\documentclass[12pt,reqno]{amsart}
\usepackage{amssymb}
\usepackage{amsmath,amscd}
\usepackage{graphicx}
\usepackage{amsthm}
\usepackage{url}
\usepackage[colorlinks=true, linkcolor=black, citecolor=black, urlcolor=blue]{hyperref}
\pdfoutput=1

\setbox0=\hbox{$+$}
\newdimen\plusheight
\plusheight=\ht0
\def\+{\;\lower\plusheight\hbox{$+$}\;}

\setbox0=\hbox{$-$}
\newdimen\minusheight
\minusheight=\ht0
\def\-{\;\lower\minusheight\hbox{$-$}\;}

\setbox0=\hbox{$\cdots$}
\newdimen\cdotsheight
\cdotsheight=\plusheight%\ht0
\def\cds{\lower\cdotsheight\hbox{$\cdots$}}

\newcommand{\cross}{\times}

\renewcommand{\(}{\left\(}
\renewcommand{\)}{\right\)}
\renewcommand{\[}{\left\[}
\renewcommand{\]}{\right\]}

\numberwithin{equation}{section}
  \theoremstyle{plain}
\newtheorem{theorem}{Theorem}[section]
\newtheorem{lemma}[theorem]{Lemma}

\newtheorem{conjecture}[theorem]{Conjecture}
\newtheorem{proposition}[theorem]{Proposition}

\newtheorem{question}[theorem]{Question}

\newtheorem{observation}[theorem]{Observation}

\newcounter{ale}
\newcommand{\abc}{\item[\alph{ale})]\stepcounter{ale}}

\newenvironment{liste}{\begin{itemize}}{\end{itemize}}
\newcommand{\aliste}{\begin{liste} \setcounter{ale}{1}}
\newcommand{\zliste}{\end{liste}}

\newenvironment{abcliste}{\aliste}{\zliste}

\begin{document}
\title{Real Places and Surface Bundles}
\author{Jonah Sinick}\thanks{The author was partially supported by the NSF}
\address{Dept.~of Mathematics, MC-382, University of Illinois, Urbana, IL 61801, USA}
\email{jsinick2@illinois.edu}
\maketitle
\begin{abstract}

Our results complement D. Calegari's result that there are no hyperbolic once-punctured torus bundles over $S^1$ with trace field having real place. We exhibit several infinite families of pairs $(-\chi, p)$ such that there exist hyperbolic surface bundles with over $S^1$ with fiber having $p$ punctures and Euler characteristic $\chi$ with trace field of having a real place. This supports our conjecture that there exist such examples for each pair $(-\chi, p)$ with finitely many known exceptions.

\end{abstract}

\section{Introduction}\label{sec:intro}

A finite volume orientable \emph{hyperbolic 3-manifold} has the form $M = \mathbb{H}^3/\Gamma$ where $\mathbb{H}^3$ is hyperbolic 3-space and $\Gamma$ is a discrete, torsion-free, finite covolume subgroup of $\text{Isom}^{+}(\mathbb{H}^3) = PSL_2(\mathbb{C})$. Such manifolds occupy a prominent role in the study of 3-manifolds in general \cite{Thu1}. By Mostow-Prasad rigidity, the topology of such an $M$ determines $\Gamma$ uniquely up to conjugation and the \emph{trace field} $\mathbb{Q}(\text{tr}(\Gamma)) = \mathbb{Q}(\{\text{tr}(\gamma): \gamma \in \Gamma\})$ is a \emph{finite} extension of $\mathbb{Q}$ which is a topological invariant of $M$. For a justification of these facts and more, see Chapter 3 of \cite{MR}. Since all hyperbolic 3-manifolds that we consider are finite volume and orientable, we now drop these modifiers.\

The trace field of a hyperbolic 3-manifold is not an abstract number field, but a concrete subfield of the complex numbers, that is, a pair $(K, \sigma)$ where $K$ is an abstract field and $\sigma: K \rightarrow \mathbb{C}$ is nonzero ring homomorphism. A short argument shows that $\sigma(K) \not \subset \mathbb{R}$ since otherwise $\Gamma$ could not have finite covolume \cite{MR}. However, it can happen that there is some other ring homomorphism $\sigma': K \rightarrow \mathbb{C}$ that $\sigma'(K) \subset \mathbb{R}$. If this is so then we say that $K$ has a \emph{real place}. The condition that $K$ has a real place is the same as the condition that a minimal polynomial for $K$ over $\mathbb{Q}$ has a real root. \

A natural question is

\begin{question}\label{neumannquestion}
Is every pair $(K, \sigma)$ with $\sigma(K) \not \subset \mathbb{R}$ the trace field of some hyperbolic 3-manifold?
\end{question}

W. Neumann has conjectured that the answer is yes, but the question is still open. An adjacent question is whether there are restrictions on the trace fields that arise from natural families of hyperbolic 3-manifolds. One of the few known results of this type is D. Calegari's intriguing theorem \cite{Cal}. Calegari's theorem is
\begin{theorem}\label{ct}
The trace field of a hyperbolic once punctured torus bundle has no real places.
\end{theorem}
Here we consider whether Calegari's result can be extended to hyperbolic surface bundles with other fibers. We show that many surfaces occur as fibers of hyperbolic surface bundles with trace field having real place.\

We now review hyperbolic surface bundles. Let $S = S(g, p)$ be the orientable, connected surface of genus $g$ with $p$ punctures. Let $\psi: S \rightarrow S$ be an orientation preserving homeomorphism. The \emph{mapping torus} of the pair $(S, \psi)$ is $M = S \times [0, 1] \slash \sim$, where $\sim$ identifies $S \cross \{0\}$ with $S \cross \{1\}$ via $\psi$. A 3-manifold is called a \emph{surface bundle over $S^1$} if it arises from this construction. If $\psi$ and $\psi'$ differ by a homeomorphism isotopic to the identity map, then the associated mapping cylinders are homeomorphic. Hence $M$ depends only on the class that $\psi$ represents in the \emph{mapping class group} $\text{Mod}(S) := \text{Homeo}^{+}(S)/\text{Homeo}_{0}(S)$; here $\text{Homeo}^{+}(S)$ is the group of orientation preserving homeomorphisms of $S$ and $\text{Homeo}_{0}(S)$ is the group of homeomorphisms of $S$ that are isotopic to the identity.\

W. Thurston showed that a mapping torus $M$ is a hyperbolic 3-manifold if and only if $\psi$ is \emph{pseudo-Anosov} \cite{Thu3}. J. Maher showed that if the Euler characteristic $\chi(S) = 2 - 2g - p < 0$ and $(g, p) \neq (0, 3)$, then a ``generic'' mapping class $\psi \in \text{Mod}(S)$ is pseudo-Anosov, so that ``most'' surface bundles over $S^1$ are hyperbolic \cite{Mah}. One motivation for studying hyperbolic surface bundles over $S^1$ is that Thurston has conjectured that every hyperbolic 3-manifold is a finite quotient of such a bundle.

In this paper we refer to the unique surface with $p$ punctures and $\chi(S) = \chi$ as the \emph{surface of type} $(-\chi, p)$, and a 3-manifold that can be realized as a surface bundle fiber of type $(-\chi, p)$ as a \emph{surface bundle of type $(-\chi, p)$}. We remark that the type of a surface bundle is not in general uniquely determined, since a given manifold can fiber over $S^1$ in multiple ways.\

In this language, Theorem \ref{ct} is that the trace field of a hyperbolic surface bundle of type $(1, 1)$ has no real places. Calegari's proof does not readily generalize to surface bundles of other types. Our goal is to substantiate
\begin{conjecture} \label{mc}
For each pair $(-\chi, p)$ with $-\chi > 1$, there exists a hyperbolic mapping torus of type $(-\chi, p)$ with trace field having a real place.
\end{conjecture}
Those surfaces with $\chi \geq 0$ and $(-\chi, p) = (1, 3)$ and are excluded since there are no hyperbolic surface bundles with these as fibers, while the surface of type $(1, 1)$ is excluded by Theorem \ref{ct}. Throughout the remainder of the paper we restrict ourselves to surfaces with $-\chi > 1$.

We were led to Conjecture \ref{mc} by utilizing existing software to collect data which led to the following.

\begin{observation}\label{ob}
For each of the following pairs $(-\chi, p)$, there exists a hyperbolic mapping torus of that type with trace field having a real place:\

\begin{center}
 {\rm
 \begin{tabular}{ r | c }
   $-\chi$ & $p$  \\ \hline
   2 & 2, 4 \\
   3 & 1, 3, 5 \\
   4 & 0, 2, 4, 6 \\   \hline
   5 & 0, 1, 3, 5, 7 \\
   6 & 0, 2, 4 \\

 \end{tabular} \hspace{1cm}
 \begin{tabular}{ r | c  }
   $-\chi$ & $p$  \\ \hline
   7 & 1, 3, 5 \\
   8 & 0, 2, 4 \\
   9 & 1, 3 \\ \hline
   10 & 0, 2, 4 \\
   11 & 1, 3 \\

 \end{tabular} \hspace{1cm}
 \begin{tabular}{ r | c  }
   $-\chi$ & $p$  \\ \hline
   12 & 0 \\
   14 & 0, 2 \\
   15 & 1 \\  \hline
   16 & 0, 2 \\
   17 & 1

 \end{tabular} \hspace{1cm}
 \begin{tabular}{ r | c  }
   $-\chi$ & $p$  \\ \hline
   18 & 0, 2 \\
   19 & 1 \\
   20 & 0, 2 \\  \hline
   21 & 1, 3 \\
   22 & 0, 2
 \end{tabular}}

\end{center}

\end{observation}

\noindent We were not able to prove our conjecture in general, but were able to prove it for several infinite families of surfaces.

\begin{theorem} \label{mt}

Conjecture \ref{mc}  is true for mapping tori of type $(-\chi, p)$ if

\begin{abcliste}
\abc $p \geq 4$ and $2|(-\chi)$
\abc $p = 0$
\abc $5p \leq -\chi$
\end{abcliste}

\end{theorem}

We use a simple fact in the proof of Theorem \ref{mt}.

\begin{proposition} \label{prop}

If $M$ is a finite volume hyperbolic 3-manifold with trace field having real place then any finite cover $M'$ of $M$ is a finite volume hyperbolic 3-manifold with trace field having real place.

\end{proposition}

The volume of $M'$ is finite since $\text{Vol}(M') = d \cdot \text{Vol}(M)$ where $d$ is the degree of the covering map $M' \rightarrow M$. Since the trace field of a finite cover of a hyperbolic 3-manifold is a subfield of the trace field of the base and a subfield of a field $K \subset \mathbb{R}$ is itself contained in $\mathbb{R}$, the truth of the Proposition \ref{prop} is immediate. So to show that Conjecture \ref{mc} is true for a surface of type $(-\chi, p)$, it suffices to show that a mapping torus of type $(-\chi, p)$ covers a hyperbolic 3-manifold with trace field having real place.\

In Section \ref{sec:twistknots} we prove a special case of Theorem \ref{mt} (a) by combining J. Hoste and P. Shanahan's explicit determination of the trace fields of hyperbolic twist knots \cite{HS} with the constructions of fibered covers of hyperbolic twist knots due to G. Walsh \cite{Wal}. In Section \ref{sec:covers} we give a criterion for a surface bundle of type $(-\chi, p)$ to be covered by a surface bundle of type $(-\chi', p')$ and use this to prove parts (a) and (b) of Theorem \ref{mt}. In Section \ref{sec:thurstonnorm} we use properties of the Thurston norm on the cohomology of a 3-manifold to prove that a particular 3-manifold with trace field having real place is a surface bundle of type $(-\chi, p)$ for many pairs $(-\chi, p)$ and this gives a proof of Theorem \ref{mt} (c). In Section \ref{sec:data} we offer empirical data substantiating Observation \ref{ob}, discuss patterns that we observed and state some open questions.

\section{Fibered Covers of Twist Knots}\label{sec:twistknots}

In \cite{HS} Hoste and Shanahan compute the trace fields of an infinite family of 3-manifolds, namely the complements of the twist knots in $S^3$. The twist knot $K_m$ is defined by Figure 1 for $m > 0$. When $m = 2$, $K_m$ is the trefoil knot and when $m = 3$, $K_m$ is the figure-eight knot. For $m \geq 3$, $S^3/K_m$ is hyperbolic so that the trace field is well-defined. Hoste and Shanahan found that if $m \geq 3$ and $m$ is even then the trace field of $S^3/K_m$ is of odd degree over $\mathbb{Q}$. As a polynomial of odd degree has a real root, the trace field of such a knot complement $S^3/K_m$ has a real place.\\

Every two-bridge knot is associated to a fraction $\frac{p}{q}$, written in lowest reduced terms and with $q > 0$. In \cite{Wal}, Walsh showed that the complement of such a knot is finitely covered by surface bundle of type $(q - 3, q - 1)$. As Figure 1 shows, $K_m$ is the two-bridge knot associated to the rational number $\cfrac{1}{2+\frac{1}{m - 1}} = \frac{m - 1}{2m -1}$. Since $\gcd(m - 1, 2m - 1) = 1$, it follows that $S^3/K_m$ is a finitely covered by a surface bundle of type $(2m - 4, 2m - 2)$.\

By Proposition \ref{prop}, the result in \cite{HS}  implies that this surface bundle has a real place if $m$ is even. This proves Conjecture \ref{mc} for surfaces of type $(4n - 4, 4n - 2)$ for $n \geq 2$. As we will see in Section \ref{sec:covers}, by utilizing another method we can prove Theorem \ref{mt} (b) which is a strictly stronger statement.\

We remark that since Leininger \cite{Lein} showed that the twist knot complements are virtually fibered by punctured torus bundles, one can use Leininger's work in conjunction with Hoste and Shanahan's work to exhibit punctured torus bundles with trace field having real place. However, though we have not checked this rigorously, it appears that following this strategy can only prove Conjecture \ref{mc} for a subset of the family of surfaces in Theorem \ref{mt} (b).

 \begin{figure}
\includegraphics[scale=0.4]{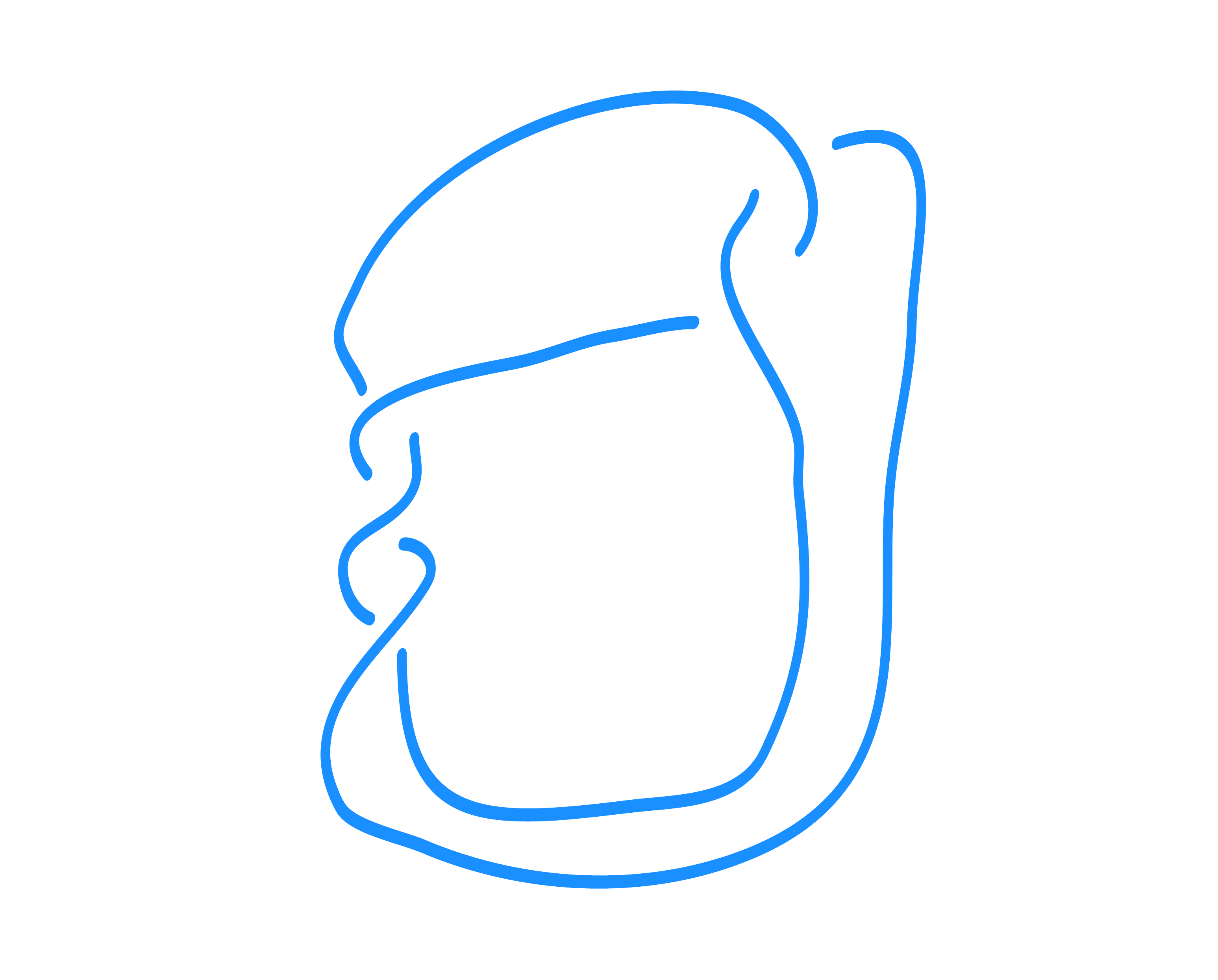}
  \caption{A depiction of $K_m$ for $m = 4$. Hoste and Shanahan's $K_m$ is isomorphic to the knot depicted with $m - 1$ crossings on the left hand side of the figure.}
\end{figure}

\section{Covers of Surface Bundles}\label{sec:covers}

In this section we prove Theorem \ref{mt} (a) and (b). First we prove

\begin{theorem} \label{bc} A surface bundle of type $(-\chi, p)$ is finitely covered by a surface bundle of type $(-\chi', p')$ if $\chi'/\chi = d \in \mathbb{N}$ and $p \leq p' \leq dp$.
\end{theorem}

After doing this we exhibit surface bundles of types $(2, 4)$ and $(2, 0)$ with trace field having real place. Parts (a) and (b) of Theorem \ref{mt} will then follow easily from Theorem \ref{bc}.\

To prove Theorem \ref{bc} we need a theorem of Massey \cite{Mas} characterizing when one punctured surface covers another.

\begin{theorem} \label{pc} For the surface of type $(-\chi', p')$ to be a degree $d$  cover of the surface of type $(-\chi, p)$ it is necessary and sufficient that (i) $\frac{-\chi'}{\chi} = d \in \mathbb{N}$ and (ii) $p \leq p' \leq pd$.
\end{theorem}

Massey's proof proceeds by cut and paste arguments. In order to use our knowledge of covers of surfaces of type $(-\chi, p)$ to pass to knowledge about covers of surface bundles of type $(-\chi, p)$ we need another lemma.

\begin{lemma} \label{mancov} Suppose $M$ is a surface bundle with fiber $S$ and that $S$ is finitely covered by $S'$. Then $M$ has a finite cover $\tilde{M}$ which is a surface bundle with fiber $S'$.
\end{lemma}

\noindent \emph{Proof of Lemma \ref{mancov}}: Let $\psi: S \rightarrow S$ define a mapping torus homeomorphic to $M$,  and let $\rho: S' \rightarrow S$ be a covering map. After picking base points for fundamental groups, we have the inclusion $\rho_{*}(\pi_1(S') \subset \pi_1(S)$. The map $\psi_{*}: \pi_1(S) \rightarrow \pi_1(S)$ is a bijection and so maps index $d$ subgroups to index $d$ subgroups. Since the number $m$ of degree $d$ subgroups of $\pi_1(S)$ is finite, $(\psi_{*})^m = (\psi^m)_{*}$ fixes $\rho_{*}(\pi_1(S')$. Define $M_1$ to be the degree $m$ cover of $M$ that is the mapping torus defined by $\psi^m$, corresponding to ``unwinding'' $M$ in the direction transverse to the fiber $m$ times. According to the lifting criterion (Prop. ~1.33 of \cite{Hat}) there exists a lift $\widetilde {\psi^m}$ of $\psi^m$ so that the following diagram commutes.

\begin{center}
$\begin{CD}
S' @>\widetilde {\psi^m}>> S\\
@VV \rho V @VV \rho V\\
S @>\psi^m>> S
\end{CD}$
\end{center}

Then the surface bundle $\widetilde {M}$ with fiber $S'$ and monodromy $\widetilde {\psi^m}$  is a cover of $M$ of degree $dm$. \qedsymbol \\

\noindent \emph{Proof of Theorem \ref{bc}}: This is immediate from Theorem \ref{pc} and Lemma \ref{mancov}. \qedsymbol \\

\noindent \emph{Proof of Theorem \ref{mt} (a)}: Let $S$ be the four-punctured sphere, that is, the surface of type $(2, 4)$. Let $S'$ be a surface of type $(-\chi', p')$, where $p' \geq 4$ and $2|(-\chi')$, as in the statement of Theorem \ref{mt} (a). Then since $\chi' = 2 - 2g' - p'$, $p' \leq -\chi' + 2$ so that by Theorem \ref{bc}, to prove Theorem \ref{mt} (a) it suffices to exhibit an $S$-bundle with trace field having a real place. An $S$-bundle can be specified by giving the mapping class $[\psi] \in \text{Mod}(S)$. The subgroup of $\text{Mod}(S)$ fixing one of the four punctures is isomorphic to the braid group on three strands. This group is generated by the braids $\sigma_1$ and $\sigma_2$ in Figure 2. We take $M$ to be the $S$-bundle given by the mapping class $[\psi]$ corresponding to the braid ${\sigma_1}{\sigma_2}^4{\sigma_1}^{-2}{\sigma_2}{\sigma_1}{\sigma_2}^{-1}{\sigma_1}{\sigma_2}$. SnapPy \cite{CDW} computes a hyperbolic structure on the resulting $S$-bundle, and Snap \cite{CGHN} rigorously proves that the trace field of $M$ is generated by the roots of the polynomial $x^4 - x^3 - 2x^2 - x + 1$. Since this polynomial has two real roots, the trace field of $M$ has a real place.\

Now Theorem \ref{mt} (a) follows from the existence of a surface bundle of type $(2, 4)$ with trace field having a real place and Theorem \ref{bc}. \qedsymbol \\

 \begin{figure}
\includegraphics[scale=0.9]{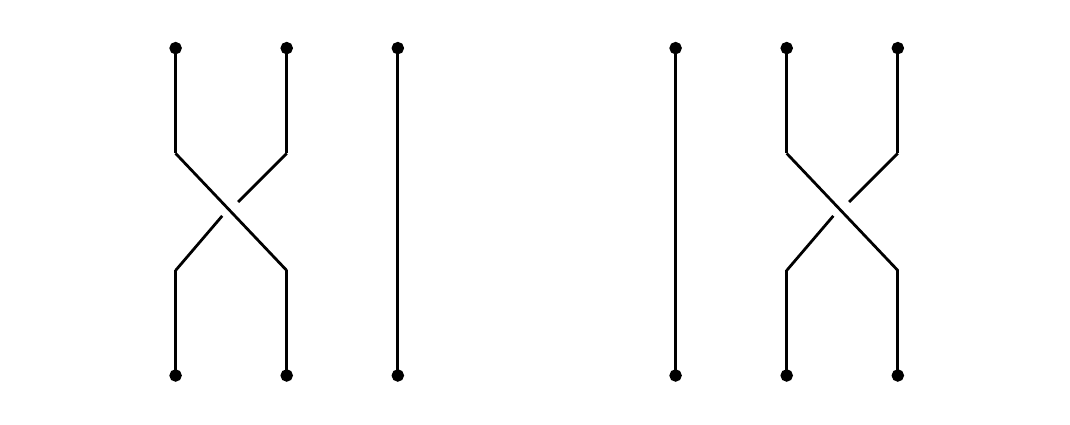}
  \caption{A depiction of the braids $\sigma_1$ (on the left) and $\sigma_2$ (on the right)}
\end{figure}

\noindent \emph{Proof of Theorem \ref{mt} (b)}: Let $S$ be the closed surface of genus $2$, that is, the surface of type $(2, 0)$. As above, we specify an $S$-bundle by specifying a mapping class in $[\psi] \in \text{Mod}(S)$. The group $\text{Mod}(S)$ is generated by Dehn twists in 5 standard curves \cite{Bir}. We order these curves according to $x$ value on pg. ~169 of \cite{Bir} as $a$, $b$, $c$, $d$ and $e$. Let $\tau_{\alpha}$ be the mapping class of the Dehn twist about the curve $\alpha$. Let $M$ be the $S$-bundle corresponding to $\tau_{a}(\tau_{b} \tau_{c} \tau_{d} \tau_{e} \tau_{b})^{-1}$. Twister \cite{Hal} generates a triangulation of $M$, and Snap \cite{CGHN} rigorously proves that the trace field of $M$ is generated by the roots of the polynomial $x^8 - 5x^6 + 12x^4 - 9x^2 + 2$. Since this polynomial has four real roots, the trace field of $M$ has a real place. \

Theorem \ref{mt} (b) follows from the existence of a surface bundle of type $(2, 0)$ with trace field having a real place and Theorem \ref{bc}. \qedsymbol \\

\section{A 3-Manifold Fibering In Many Ways}\label{sec:thurstonnorm}

In this section we prove

\begin{theorem} \label{wf}

There is a hyperbolic 3-manifold with trace field having real place that is a surface bundle of type $(-\chi, p)$ for every pair with $5p \leq -\chi$ coming from a surface except for $(5, 1)$.
\end{theorem}

Since according to Observation \ref{ob}, Conjecture \ref{mc} is true for $(-\chi, p) = (5, 1)$; Theorem \ref{wf} implies Theorem \ref{mt} (c).\

If a hyperbolic 3-manifold $M$ fibers over $S^1$, then the 1-form $\frac{1}{2\pi} d \theta$ on $S^1$ pulls back to a cohomology class $\phi \in H^1(M; \mathbb{Z})$, and we say that $\phi$ induces the fibration. In \cite{Thu2}, Thurston provides an elegant characterization of all fibrations of $M$ over $S^1$. First Thurston defines a norm $||*||_T$ (now called the \emph{Thurston norm}) on  $H^1(M; \mathbb{R})$. We do not need the definition of the Thurston norm here but will use the following consequence of the definition.

\begin{proposition} \label{prop2}

If $\phi$ induces a fibration of a hyperbolic 3-manifold $M$ over $S^1$ with fiber $S$ then $||\phi||_T = -\chi$ where $-\chi$ is the Euler characteristic of $S$.

\end{proposition}

In \cite{Thu2} Thurston shows that the unit norm ball of $||*||_T$ is a convex polyhedron symmetric about the origin and with coordinates of the vertices in $\mathbb{Q}$.  We denote its unit norm ball by $B_T$. Setting $b_1(M) = \text{dim}(H^{1}(M; \mathbb{R}))$, we say that a face of $B_T$ is \emph{top dimensional} if its dimension is $b_1(M) - 1$. For $\phi \in H^1(M; \mathbb{R})$, the ray from the origin to $\phi$ intersects the interior of a unique face of $B_T$, and we say that $\phi$ \emph{lies above} that face. A face and its reflection through the origin play the same role, so we abuse notation, referring to a \emph{pair} of opposite faces of the unit norm ball as a face. The main theorem of \cite{Thu2} is:

 \begin{theorem}\label{tt}
 There is a collection of top dimensional faces of $B_T$ called the fibered faces so that $\phi \in H^1(M; \mathbb{Z})$ induces a fibration of $M$ if and only if $\phi$ lies over one of the faces in this collection.
 \end{theorem}

 Thurston's theorem is trivial for $M$ if $b_1(M) = 1$, but if $b_1(M) > 1$, it shows that if $M$ is a surface bundle then there are infinitely many pairs $(-\chi, p)$ such that $M$ is a surface bundle of type $(-\chi, p)$. In this section we give an example of a hyperbolic 3-manifold $M$ with trace field having real place and $b_1(M) > 1$, and work out the pairs $(-\chi, p)$ such that $M$ is a surface bundle of type $(-\chi, p)$.\

For our manifold, we take $M$ as in Figure 3. We found this manifold in Morwen Thistlethwaite's recent census of cusped hyperbolic 3-manifolds that are gluings of $8$ ideal tetrahedra. Thistlewaite's data has been incorporated into the most recent version of Snappy \cite{CDW}.

 \begin{figure}
\includegraphics[scale=0.4]{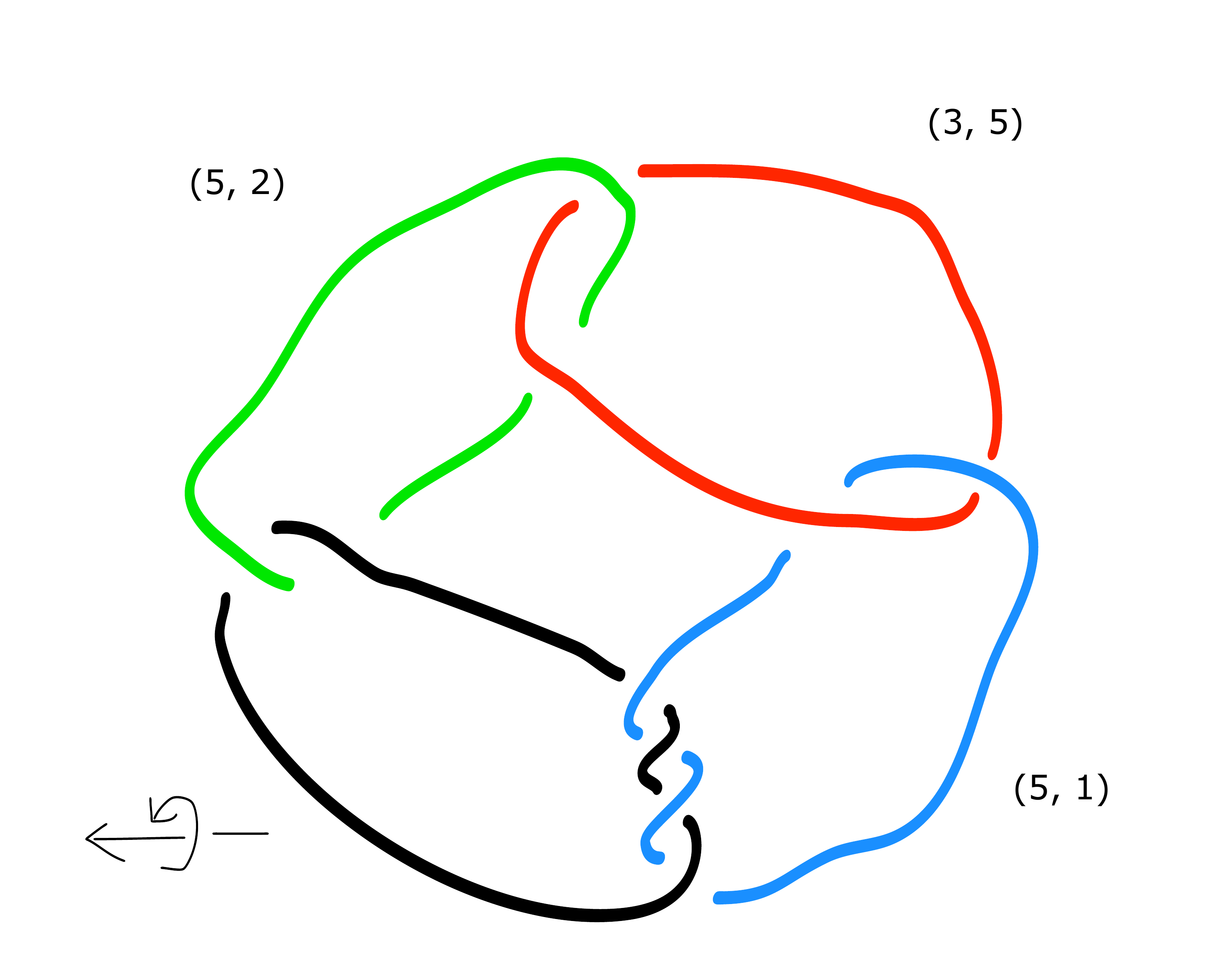}
  \caption{The ordered pairs adjacent to each of the colored loops denote Dehn surgery coefficients, with Dehn surgery convention as indicated in the lower lefthand corner.}
\end{figure}

Snap \cite{CGHN} computes that a minimal polynomial for the trace field of $M$ is $x^6 - 2x^4 + 4x^2 - 1$. Since this polynomial has two real roots, the trace field of $M$ has a real place. To study the surfaces $S$ such that $M$ is a surface bundle with fiber $S$ we need to compute $B_T$ and determine which faces of $B_T$. This will allow us to use Proposition \ref{prop2} to determine the Euler characteristics of the fibers. First we show:

\begin{theorem}\label{ff}\ Let $B_T$ be the unit ball of the Thurston norm on $H^1(M; \mathbb{R})$. Then all top dimensional faces of $B_T$ are fibered.\
\end{theorem}

\noindent \emph{Proof of Theorem \ref{ff}}: Stallings cite{St} gave an algebraic criterion for a class $\phi \in H^1(M; \mathbb{Z})$ to induce a fibration of $M$ over $S^1$. Under the assumption that $\pi_1(M)$ is specified by a presentation with two-generators and one-relator presentation, Brown \cite{Br} gave an algorithm for determining all $\phi \in H^1(M; \mathbb{Z})$ that satisfy the hypothesis of Stallings' criterion. We use Brown's algorithm to apply Stallings' criterion and determine all $\phi \in H^1(M; \mathbb{Z})$ that induce fibrations of $M$. For a nice exposition of the perspective on these topics that we use here, see sections 4 and 5 of \cite{DT}.\

 SnapPy \cite{CDW} computes a presentation of $\pi_1(M)$ given by
\begin{align*}
 <\alpha, \beta \hspace{0.15 cm}|\hspace{0.15 cm} {\alpha}^2 {\beta}^3 {\alpha}^2 {\beta}^{-2}{\alpha}^{-3} {\beta}^{-2} {\alpha}^2 {\beta}^3 > .
\end{align*}

 \begin{figure}
\includegraphics[scale=0.7]{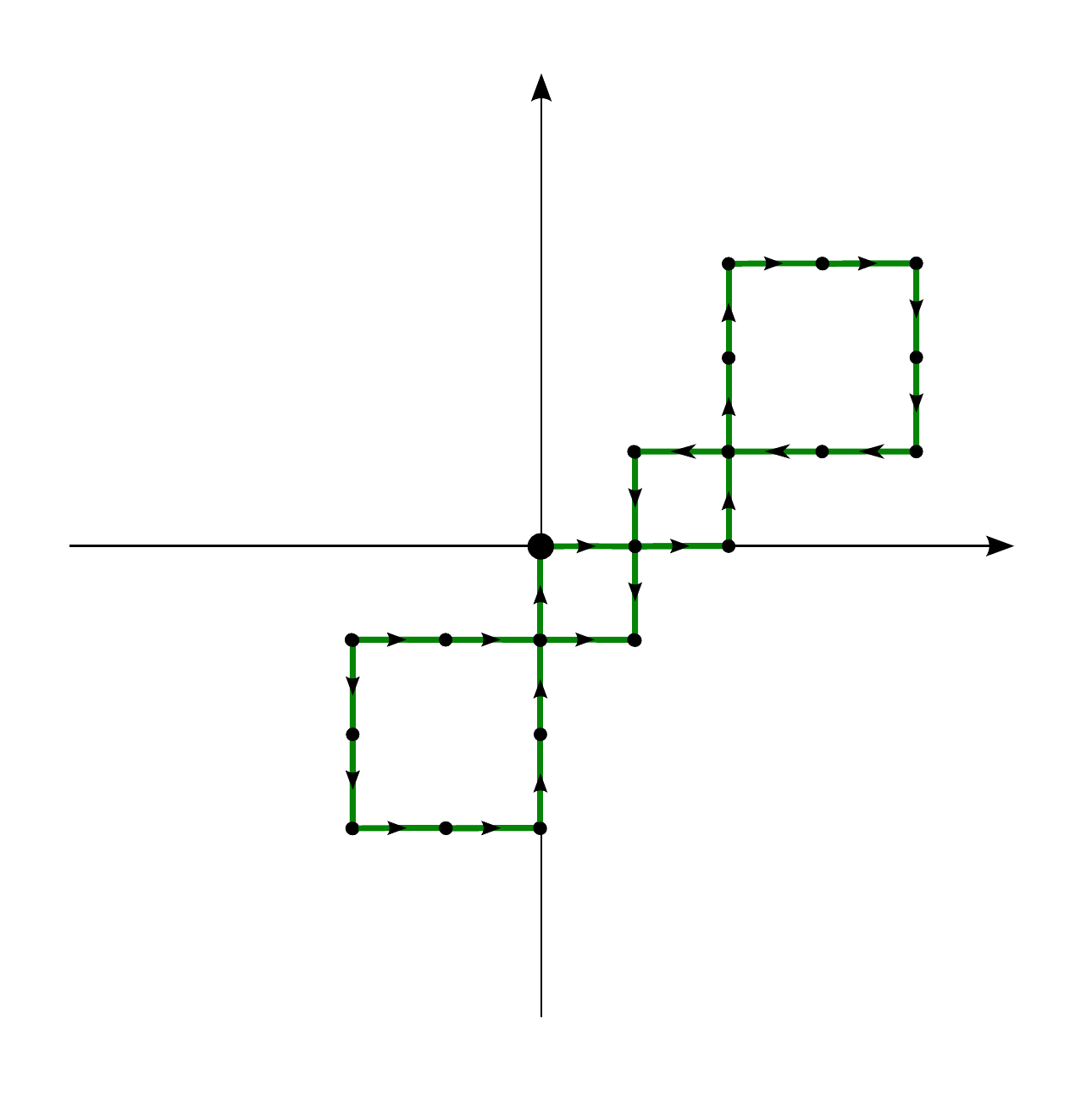}
  \caption{Path associated to presentation of $\pi_1(M)$. Defined only up to translation by a lattice point.}
\end{figure}

\noindent Since this presentation is 2-generator and 1-relator, Brown's criterion can be applied. Brown's criterion is stated in terms of a bounded path $P$ in $\mathbb{R}^2$ consisting of horizontal and vertical segments that can be produced from a presentation of a 2-generator 1-relation group $G$. The path starts at the origin and is formed by reading the relator from left to right, drawing a horizontal unit segment for each appearance of $\alpha$ and a vertical unit segment for each appearance of $\beta$. We give a depiction of the path associated to the above presentation of $\pi_1(M)$ in Figure 4.

Brown's criterion is that $\phi = (c, d) \in H^1(M; \mathbb{Z})$ induces a fibration of $M$ if and only if as $e$ varies over elements of $\mathbb{R}$, the largest and smallest value of $e$ for which the line $L_e$ given by $cx + dy = e$ intersects $P$ each give a line which intersects $P$ exactly once. Here $L_e$ is considered to intersect $P$ multiple times if $P$ crosses over itself where it meets $L_e$, but $L_e$ is considered to intersect $P$ only once if $L_e$ is horizonal (resp. ~vertical) and $L_e$ intersects a single horizontal (resp. ~vertical) segment. Applying Brown's criterion to Figure 4, it is clear that the only pairs $(c, d) \in H^1(M; \mathbb{Z})$ that do not induce fibrations of $M$ are those lie on the three lines that go through the origin and have slope $0$, $-1$ and $\infty$. By Theorem \ref{tt}, there cannot be any top dimensional faces of $B_T$ that are not fibered, completing the proof of Theorem \ref{ff}. \qedsymbol\

We can use Theorem \ref{ff} to determine the precise shape of $B_T$.

\begin{theorem}\label{tnb}\ The unit ball of the Thurston norm $B_T$ is as pictured in Figure 5.
\end{theorem}

 \begin{figure}
\includegraphics[scale=0.7]{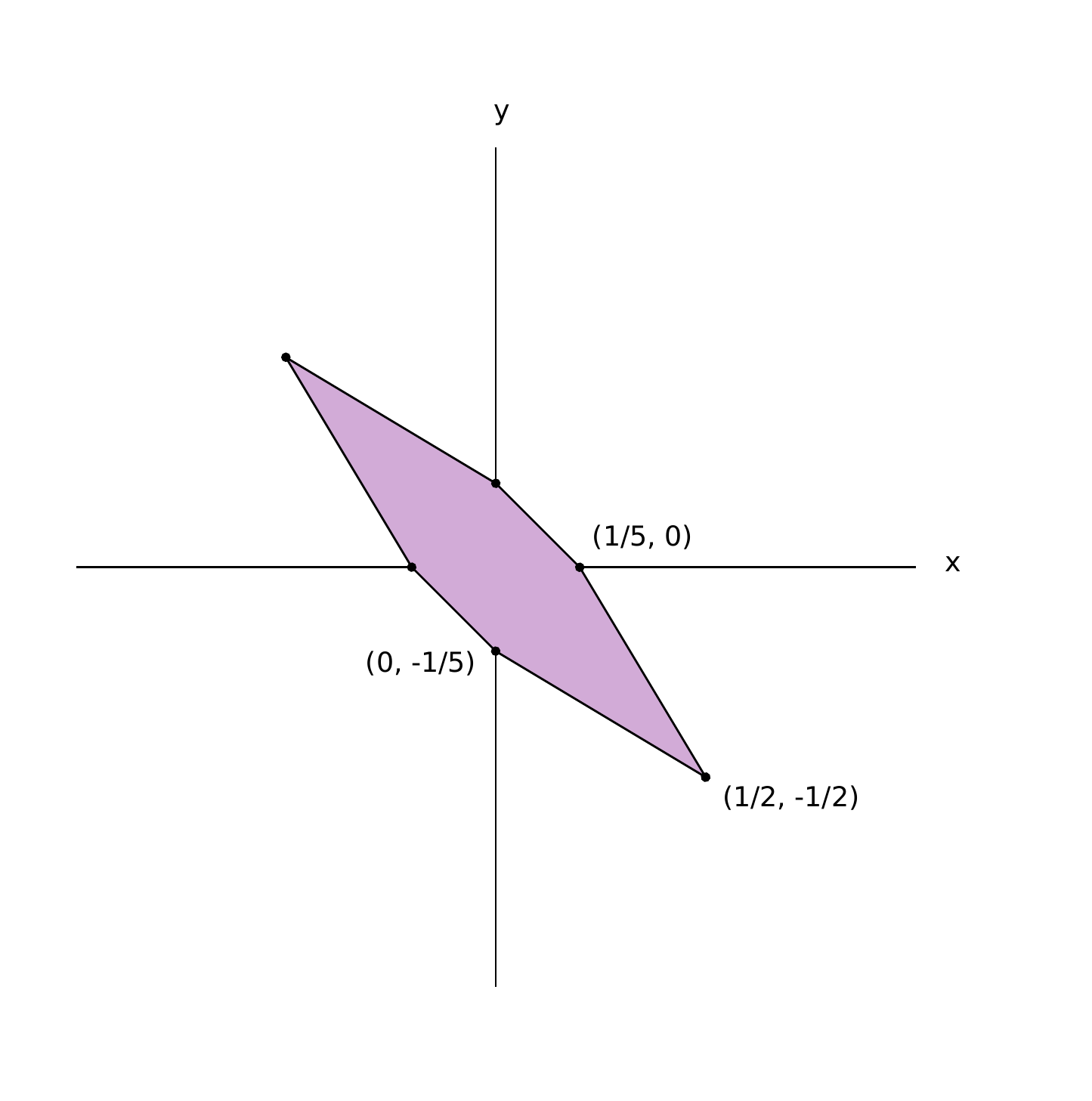}
  \caption{The unit ball $B_T$ of the Thurston norm on $H^1(M; \mathbb{R})$}
\end{figure}
\noindent \emph{Proof of Theorem \ref{tnb}}: We use Theorem \ref{ff} together with properties of the Alexander norm on $H^1(M; \mathbb{R})$. The Alexander norm $||*||_A$ on  $H^1(M; \mathbb{R})$ was introduced by McMullen in \cite{McM}. McMullen showed that if $N$ is a 3-manifold and $\phi \in H^1(N; \mathbb{Z})$ then $||\phi||_A \leq ||\phi||_T$ where $||*||_T$ is the Thurston norm, and that equality holds if $\text{dim}(H^1(M; \mathbb{R})) > 1$ and $\phi$ induces a fibration of $N$. Applying this result to $M$, $||\phi||_A = ||\phi||_T$ for each $\phi$ not lying on the lines through the origin with slope $0$, $-1$ and $\infty$. Since McMullen \cite{McM} and Thurston \cite{Thu2} showed that their respective norms are linear on rays through the origin and continuous on $H^1(M; \mathbb{R})$, we must have $||\phi||_A = ||\phi||_T$ identically. So $B_T$ = $B_A$ where $B_A$ is the unit ball of the Alexander norm on $H^1(M; \mathbb{R})$.\

The advantage of working with the Alexander norm is that it is easy to compute. Indeed, the unit norm ball on $H^1(M; \mathbb{R})$ is defined to be the dual of the Newton polytope of the Alexander polynomial of $N$. The Alexander polynomial of $M$ is
\begin{align*}
a^5 b^5 + a^5 b^4 + a^4 b^5 + a^4 b^4 - a^3 b^3 + a^2 b^2 - ab - a - b - 1
\end{align*}

Thus, its Newton polytope has vertices $(0, 0)$, $(1, 0 )$, $(0, 1)$, $(4, 5)$, $(5, 4)$ and $(5, 5)$. The dual of this polytope has vertices given by \noindent$(0, 1/5)$, $(1/5, 0)$, $(-1/2, 1/2)$, $(-1/5, 0)$, $(0, -1/5)$, and  $(1/2, -1/2)$ (Figure 5). This completes the proof of Theorem \ref{tnb} . \qedsymbol\\

Using our computation of $B_T$, the fact that the Thurston norm is linear on rays, and Proposition \ref{prop2}, if $\phi$ induces a fibration of $M$ then we can compute the Euler characteristic of the corresponding fiber. To completely determine the fiber we need to be able to compute the number of punctures of the fiber. To this end we prove:

\begin{theorem}\label{np}\ If $\psi = (x, y)$ lies over a fibered face of $B_T$ then $\psi$ induces fibration of $M$ of type $(-\chi, p)$, where $p$ = $|x + y|$.
\end{theorem}

\noindent  \emph{Proof of Theorem \ref{np}}: Since $M$ is a one-cusped hyperbolic 3-manifold, removing a neighborhood of the cusp of $M$ gives a $3$-manifold $M'$ with boundary $\partial M'$ consisting of a torus. Pick $\phi$ inducing a fibration of $M$. Pick generators $\mu$ and $\lambda$ of $H_1(\partial M'; \mathbb{Z})$  with intersection number $1$ such that $i(\mu)$ is transverse to the fiber and $i(\lambda)$ is parallel to the fiber in the fibration of $M$ induced by $\phi$. Then $\phi(\mu) = p$ and $\phi(\lambda) = 0$, where $p$ is the number of punctures of a fiber. In some sense this gives a formula for $p$, but this formula depends on a choice of generators of $H_1(\partial M'; \mathbb{Z})$ which is special to $\phi$. A formula for $p$ that does not depend on generators of $\phi$ is given by
\begin{align}
|\gcd(\phi(\mu), \phi(\lambda))| = p.
\end{align}
To see that (4.1) is true, first note that it is true if $\mu$ and $\lambda$ are generators that depend on $\phi$ as above. Now, any choice of generators of $H_1(\partial M'; \mathbb{Z}) \cong \mathbb{Z}^2$ differs from that of $\mu, \lambda$ by an invertible linear change of coordinates, that is, an element of $GL(2, \mathbb{Z})$. Since the $\gcd$ of two numbers is invariant under $
\left[ \begin{array}{cc}
1 & 1 \\
0 & 1 \\
\end{array} \right]$,  $\left[ \begin{array}{cc}
0 & -1 \\
1 & 0 \\
\end{array} \right]$, and
   $\left[ \begin{array}{cc}
1 & 0 \\
0 & -1 \\
\end{array} \right]$
 and these matrices generate $GL(2, \mathbb{Z})$, equation (4.1) holds for any choice of generators of $H_1(\partial M; \mathbb{Z})$.

With its preferred choice of generators of $H_1(\partial M'; \mathbb{Z})$, SnapPy \cite{CDW} computes that their images in $H_1(M; \mathbb{Z})/(\text{torsion})$ to be $\mu = (-4, -4)$ and $\lambda = (-5, -5)$. Writing $\phi = (x,y)$, we have $\gcd(\phi(\mu), \phi(\lambda)) = \gcd(-4x - 4y, -5x - 5y) = x + y$. Thus Theorem \ref{np} follows from (4.1). \qedsymbol\\

\noindent \emph{Proof of Theorem \ref{wf}}: We fix $p > 0$ and determine those $\chi$ such that $M$ is a surface bundle of type $(-\chi, p)$. By Theorem \ref{np}, such $\chi$ come from classes $(x, y)$ with $x + y = \pm p$. Since the classes $(x, y)$ and $(-x, -y)$ are the same up to orientation, there is no loss of generality in restricting ourselves to those classes lying on the line $x + y = p$. Since the Thurston norm is symmetric about the lines $y = x$ and $x = 0$ and those classes lying on these lines do not induce fibrations of $M$, we can and do further restrict our attention to those classes $(x, y)$ with $x > 0$ and $y >  x$. Of these classes, those on the line $y = 0$ do not induce fibrations of $M$, but all others do. Recall that when a class induces a fibration, the Thurston norm computes its Euler characteristic. The formula for the Thurston norm of $(x, y)$ depends on whether $ y < 0$ or $y > 0$. Thus we divide our analysis into two cases.\\

\noindent \emph{Case 1}: Suppose $y > 0$. Then $(x, y)$ lies over the face of $B_T$ determined by the two vertices $(1/5, 0)$ and $(0, 1/5)$, so in this case
\begin{align*}
||(x, y)||_T = 5x + 5y = x + 5(p - x) = 5p.
\end{align*}
\noindent So every class $(x, y)$ with $y > 0$ such that $y < x$ and $x + y = p$ gives a realization of $M$ as a surface bundle of type $(5p, p)$. The pairs $(-\chi, p)$ such that there is a class $(x, y)$ with $y > 0$ that induces a fibration of $M$ of type $(-\chi, p)$ are those pairs $(5p, p)$ for $p \geq 2$.\\

\noindent \emph{Case 2}: Suppose $y < 0$.  Then $(x, y)$ lies in the cone over the face of the Thurston norm ball determined by $(1/5, 0)$ and $(1/2, -1/2)$, so in this case
\begin{align*}
||(x, y)||_T = 5x + 3y  = 5x + 3(p - x) = 2x + 3p.
\end{align*}
\noindent In this case the condition $y < x$ is automatically satisfied since $x + y = p$. The condition $y < 0$ is equivalent to $p < x$ which is in turn equivalent to
\begin{align*}
5p < 2x + 3p = ||(x, y)||_T
\end{align*}
\noindent Recall that for any pair $(-\chi, p)$ that comes from a surface, $p$ and $\chi$ have the same parity. So for any pair $(-\chi, p)$ with $5p < -\chi$ we can find a class $(x, y)$ with $y < 0$ that induces a fibration of $M$ of type $(-\chi, p)$.\\

From our above analysis we see that $M$ is a surface bundle of type $(-\chi, p)$ precisely when $5p  \leq  \chi$ and $(-\chi, p) \neq (5, 1)$ as claimed. \qedsymbol

\section{Empirical Data and Open Questions}\label{sec:data}

Here we give a table of cusped census manifolds with trace fields having real places which are surface bundles of type $(-\chi, p)$ for various pairs $(-\chi, p)$.  This table together with Theorem \ref{mt} (a) and (b) substantiates Observation \ref{ob}. We compiled it using N. Dunfield's data \cite{Dun} and Snap \cite{CGHN}. The manifolds are labelled by their labels in the cusped census \cite{CHW}.\\

\vspace{0.5 cm}

 {\rm
 \begin{tabular}{ r | c }
   $(-\chi,p)$ & Manifold  \\ \hline
   (2, 2) & v2640 \\
   (3, 1) & m036 \\
   (3, 3) & s493 \\ \hline
   (3, 5) & m043 \\
   (4, 2) & v2677 \\
   (4, 4) & s500 \\

 \end{tabular} \hspace{1cm}
 \begin{tabular}{ r | c  }
   $(-\chi,p)$ & Manifold  \\ \hline
   (4, 6) & s147\\
   (5, 1) & m034 \\
   (5, 3) & v0163 \\ \hline
   (5, 5) & v0003 \\
   (5, 7) & m172 \\
   (6, 2) & v3212\\

 \end{tabular} \hspace{1cm}
 \begin{tabular}{ r | c  }
   $(-\chi,p)$ & Manifold  \\ \hline
   (6, 4) & s500 \\
   (7, 1) & m078 \\
   (7, 3) & v3238 \\ \hline
   (7, 5) & v0022 \\
   (8, 2) & m297 \\
   (8, 4) & v3193 \\

 \end{tabular} \hspace{1cm}

\vspace{1 cm}

 \begin{tabular}{ r | c  }
   $(-\chi,p)$ & Manifold \\ \hline
   (9, 1) & m011 \\
   (9, 3) & v0170\\
   (10, 2) & m200\\ \hline
   (10, 4) & v1251 \\
   (11, 1) & m019\\
   (11, 3) & v1721\\

 \end{tabular} \hspace{1cm}
 \begin{tabular}{ r | c  }
   $(-\chi,p)$ & Manifold  \\ \hline
   (14, 2) & s156\\
   (15, 1) & m070 \\
   (16, 2) & s550\\ \hline
   (17, 1) &  m044\\
   (18, 2) & s133\\

  \end{tabular} \hspace{1cm}
 \begin{tabular}{ r | c  }
   $(-\chi,p)$ & Manifold  \\ \hline
   (19, 1) & m055\\
   (20, 2) & s457 \\ \hline
   (21, 1) & m064\\
   (21, 3) & v1609\\
   (22, 2) & v2603 \\
 \end{tabular}\\

\vspace{0.5 cm}

While investigating the trace fields of the fibered manifolds in cusped census, we found that the trace fields attached to hyperbolic surface bundles of type $(-\chi, p)$  for $(-\chi, p) \neq (1, 1)$ do not appear to exhibit any regularities. We are thus led to ask

\begin{question}\label{newneumannquestion} Suppose $(-\chi, p) \neq (1, 1)$. Is every pair $(K, \sigma)$ with $\sigma(K) \not \subset \mathbb{R}$ the trace field of some hyperbolic surface bundle of type $(-\chi, p)$?
\end{question}

Empirically, it seems that the trace field of a small volume cusped hyperbolic surface bundle typically has about $1/3$ as many real places as its degree over $\mathbb{Q}$. Since a complex place corresponds to a conjugate pair of nonreal roots of a polynomial, this observation suggests the possibility that a ``random'' place of a ``random'' cusped hyperbolic 3-manifold is real with probability $1/2$. We therefore ask the following pair of questions:

\begin{question}\label{RealPlaceQuestion} Is there a natural notion of ``random number field'' such that the expected fraction of places of a random number field that are real is $1/2$?
\end{question}

\begin{question}\label{RealPlaceQuestion} Assuming that a notion of ``random number field'' as in the preceding question is found, can one find a natural notion of ``random cusped hyperbolic surface bundle'' such that the trace fields of such manifolds are modeled well by this notion of ``random number field'' (with the restriction that such number fields have at least one complex place) and deduce that the fraction of places of a``random cusped hyperbolic surface bundle'' that are real is $1/2$?
\end{question}

\section{Acknowledgements}

The author acknowledges support from National Science Foundation grant DMS 08-38434 ``EMSW21-MCTP: Research Experience for Graduate Students.'' He thanks Nathan Dunfield for many helpful conversations and suggestions. He thanks Chris Leininger for teaching him about the Thurston norm.

\end{document}